\documentclass[12pt]{iopart}
\usepackage{iopams}  
\usepackage{amsthm}
\usepackage{bm}
\usepackage{amsfonts}
\usepackage{xcolor}
\usepackage[pdftex]{graphicx}
\usepackage{blkarray}
\usepackage{hyperref}
\usepackage[hypcap]{caption}
\newcommand{\matindex}[1]{\mbox{\scriptsize#1}}
\begin{document}

\title[Non-symmetric finite networks: the two-point resistance]{Non-symmetric finite networks: the two-point resistance}

\author{Viera \v Cerna\v nov\' a$^1$ and Juraj Brenku\v s$^2$}

\address{$^1$ Institute of Computer Science and Mathematics, Slovak University of Technology, Bratislava, Slovak Republic}
\address{$^2$ Institute of Electronics and Photonics, Slovak University of Technology, Bratislava, Slovak Republic}
\ead{\{vieracernanova\}@hotmail.com}

\begin{abstract}
\\
An explicit formula for the resistance between two nodes in a network with a non-symmetric Laplacian matrix {\bfseries L} is obtained. 
This is of great advantage e.g. in electronic circuit fault analysis, where non-linear systems have to be solved repeatedly. 
Analysis time can be greatly reduced by utilization of the obtained formula. 
The presented approach is based on the \textquotedblleft mutual orthogonality\textquotedblright{} of the full system of left and right-hand eigenvectors of a diagonalizable matrix {\bfseries L}. 
Some examples are given that demonstrate the accuracy of the approach on circuit networks.
\end{abstract}

\pacs{ 01.55.+b, 02.10.Yn, 85.40.Bh}
\ams{ 15A18, 94C12}

\maketitle

\section{Introduction}

The computation of the resistance between two nodes in a resistor network is a classical problem in
electric circuit theory, indeed it was studied by numerous authors over many years. Since electrical
potentials on a grid are governed by the same difference equations as those occurring in other
problems, the computation of the resistances can be applied to a wide range of problems, such as
lattice Green's functions, random walks, the theory of harmonic functions and others. Due to variety
of applications of the resistance problem and various methods of solution, 
this problem is still actual and intensively studied~\cite{Qwaidat2010}.

The resistance problem was often transformed in that of solving the difference equations, which is
most conveniently carried out for infinite networks. Little attention was paid to finite networks, though
these are occurring in real life. In 2004, an important result for finite lattices was published in~\cite{Wu2004}.
He expressed the resistance between two arbitrary nodes
in a resistor network in terms of eigenvalues and eigenvectors of the Laplacian matrix {\bf L} associated
with the network. In~\cite{Wu2004}, the Laplacian {\bf L} is a real symmetric $N \times N$ matrix with
$N$ real eigenvalues $\lambda_{i}$ and an orthonormal basis of eigenvectors ${\bm \psi}_{i}=(\psi_{i1},\dots,\psi_{iN})^*$. The resistance
between the nodes $\alpha$ and $\beta$ is given by
\begin{equation}
R_{\alpha\beta}=\sum_{i=2}^{N}{{1 \over \lambda_{i}}\left|\psi_{i\alpha}-\psi_{i\beta}\right|^{2}}
\label{eq:stary_wu}
\end{equation}

However, if the Laplacian matrix of the network is non-symmetric, the relation~\eref{eq:stary_wu} fails.

In this work, we study resistance networks with a real Laplacian matrix {\bf L}, symmetric or not. The
Laplacian {\bf L} is a real diagonalizable $N \times N$ matrix with $N$ real eigenvalues $\lambda_{i}$ and two basis of
\textquotedblleft mutually orthogonal\textquotedblright{} left and right-hand eigenvectors $\bm \varphi_{i}^{*}=(\varphi_{i1},\dots,\varphi_{iN})\mbox{ and }\bm\psi_{i}=(\psi_{i1},\dots,\psi_{iN})^{*}$. 
To our best knowledge, this is the first work that presents an explicit formula for the total resistance
between two nodes in a resistor network with a real Laplacian matrix, in general non-symmetric. The
obtained formula

\begin{equation}
R_{\alpha\beta}=\sum_{i=2}^{N}{{1 \over \lambda_{i}\bm\varphi_{i}^{*}\bm\psi_{i}}(\varphi_{i\alpha}-\varphi_{i\beta})(\psi_{i\alpha}-\psi_{i\beta})}
\label{eq:novy_vv}
\end{equation}
is an extension of~\eref{eq:stary_wu}. The method used in this paper is similar to that of~\cite{Wu2004}, where the orthonormal basis of eigenvectors of the matrix {\bf L} is
replaced by two \textquotedblleft mutually orthogonal\textquotedblright{} basis of left and right-hand eigenvectors of {\bf L}.

The organization of the paper is as follows. 
In section~\ref{sec:motivation} our motivation coming from the fault analysis of integrated circuits is presented and the nodal conductance matrix of a circuit, 
which is the Laplacian {\bf L}, and some of its properties are described. 
Section~\ref{sec:theoretical} gives a theoretical framework, which is based on the theory of eigenvalues and eigenvectors of real matrices. 
The existence of \textquotedblleft mutually orthogonal\textquotedblright{} basis $\bm\varphi_{1}^{*},\dots,\bm\varphi_{N}^{*}$ and $\bm\psi_{1},\dots,\bm\psi_{N}$
 of left and right-hand eigenvectors of a real diagonalizable matrix {\bf A} is proved. 
Keeping the vectors $\bm\varphi_{1}^{*},\dots,\bm\varphi_{N}^{*}$ intact, a convenient basis of right eigenvectors is constructed. 
Then, the main result is formulated as theorem and proved in the second part of section~\ref{sec:theoretical}. 
Section~\ref{sec:examples} is devoted to calculating the two-point resistance for several circuit networks.
The formula proved in section~\ref{sec:theoretical} is demonstrated on a simple bridge for two cases: simple eigenvalues and multiple eigenvalues.
Also, an example of resistance calculation for a real life circuit, an operational amplifier, is given.

\section{Fault analysis - motivation}
\label{sec:motivation}

Today's integrated circuits (ICs) represent complex systems realized on a small area using advanced fabrication technologies. At the early age of monolithic integration (early sixties), the fabrication process was dominated by yield as low as 10\% due to poor process control~\cite{Kilby2000}. 
As the fabrication process has evolved, the process control reached acceptable limits and the yield improved. On the other hand, with the increasing resolution, the fabrication process became more sensitive to impurities of the environment and had to be moved to controlled areas, 
where the ambient impurity level was under control.

These days, a typical mixed-signal IC can consist of thousands to millions of transistors, occupy an area of several tens of square mm, while the fabrication process yield with resolution as high as tens of nanometers is well above 90\%. 
Anyhow, even in a well controlled environment, the ambient still exhibits certain level of impurities. 
Besides, the fabrication process control is becoming even more complex and thus, resulting in a necessity to test the produced ICs. 

During the test development process, so called {\it fault simulations} have to be conducted to determine the fault coverage of a given test. 
If performed in a traditional way, this step might be rather time consuming, since it requires a huge amount of simulations to be executed 
due to large amount of circuit nodes and possible faults taken into account. 
During the simulation, a set of nonlinear equations describing the circuit has to be solved using an iterative algorithm, e.g. Newton-Raphson.
Application of such an algorithm to circuit analysis is described in~\cite{Vlach}.

The numerical approach to DC fault simulations, presented in this paper, evaluates the resistance between two arbitrary nodes of the circuit in a more time efficient way. This evaluation requires only one simulation for the given circuit configuration, providing information on the operating point, to be executed. It means that a system of nonlinear equations describing the circuit has to be solved only once for the given operating point. 
From this operating point information, the circuit's nodal conductance matrix can be created and the resistance between the selected two nodes can be calculated. 
As has been already shown for symmetrical matrices, such an approach can lead to a significant analysis time reduction~\cite{brenki2013}.
This resistance calculation can be used in the development of resistance-based tests, e.g. $I_{DDQ}$ testing, 
where the resistance is indirectly evaluated by measuring the current flowing from the supply source~\cite{Tsiatouhas2002, Beresinski2008}. 
For a circuit network, where active elements such as MOS transistors are present, 
the nodal conductance matrix can become non-symmetrical, since voltage-controlled current sources (VCCS) are used to model such devices, as described below.

\subsection{Nodal conductance matrix}

Formulation of the nodal conductance matrix, which is the Laplacian, is based on the Kirchhoff current law (KCL) which states the following~\cite{Kirchhoff1847}

\begin{equation}
\sum_{j=1}^{N}{}^{'} { c_{ij}(V_{i}-V_{j}) }=I_{i}\qquad i=1,2,\dots,N.
\label{eq:KCL1}
\end{equation}
Here, $c_{ij}$ is the conductance of a resistor connected to nodes $i$ and $j$, $V_{i}$ and $V_{j}$ are the electrical potentials at nodes $i$ and $j$, 
and $I_{i}$ is the current flowing into the network at node $i$. The prime denotes the omission of the term $i=j$. 
As a consequence, the following constraint applies:

\begin{equation}
\sum_{i=1}^{N} { I_{i} }=0.
\label{eq:KCL2}
\end{equation}
In general, a circuit of $N$ nodes can be described by a matrix equation of the form:

\begin{equation}
{\bf L}\vec{V}=\vec{I}
\label{eq:matrix1}
\end{equation}
or

\begin{equation}
  \begin{blockarray}{cccc}
  & & & \\
    \begin{block}{(rrrr)}
c_{11} & c_{12} & \dots & c_{1N} \\
c_{21} & c_{22} & \dots & c_{2N} \\
\vdots & \vdots & \ddots & \vdots \\
c_{N1} & c_{N2} & \dots & c_{NN} \\
    \end{block}
  \end{blockarray}
\quad
  \begin{blockarray}{c}
   \\
    \begin{block}{(c)}
V_{1} \\
V_{2} \\
\vdots \\
V_{N} \\
    \end{block}
  \end{blockarray}
  =\begin{blockarray}{c}
  \\
    \begin{block}{(c)}
I_{1} \\
I_{2} \\
\vdots \\
I_{N} \\
    \end{block}
  \end{blockarray},
\end{equation}
where ${\bf L}$ is the Laplacian (nodal conductance matrix) with elements $c_{11},\dots,c_{NN}$, $\vec V$ is a column vector of node potentials, 
and $\vec I$ is a column vector of independent current sources~\cite{Vlach}.

To compute the resistance between two arbitrary nodes $\alpha \mbox{ and } \beta$ of the network, an external source has to be connected to these nodes. Let the current flowing from the source be $I$  and the potentials at the two nodes $V_{\alpha} \mbox{ and } V_{\beta}$. Then the desired resistance is determined by the Ohm's law

\begin{equation}
R_{\alpha\beta}={V_{\alpha} - V_{\beta} \over I}.
\label{eq:Wu-6}
\end{equation}
Now, the computation is reduced to solving~\eref{eq:matrix1} for $V_{\alpha} \mbox{ and } V_{\beta}$ with the current given by

\begin{equation}
I_{i}=I(\delta_{i\alpha}-\delta_{i\beta}).
\label{eq:Wu-7}
\end{equation}
 An important property of the ${\bf L}$ matrix is 

\begin{equation}
V_{j}\sum_{i=1}^{N}{c_{ij}}=\sum_{i=1}^{N}{I_{i}}.
\label{eq:col0}
\end{equation}
From this equation follows, that since the KCL requires the sum on the right side to be zero, the sum of conductances in each column of ${\bf L}$ is equal to zero.
Another important property is given by~\eref{eq:row0}:

\begin{equation}
V_{j}\sum_{i=1}^{N}{c_{ij}}=0.
\label{eq:row0}
\end{equation}
Since $V_{j}$ is a nonzero voltage source, the sum must be equal to zero. 
Therefore, the sum of conductances in each row of ${\bf L}$ is equal to zero.

\begin{figure}[h!]
\includegraphics[width=5cm]{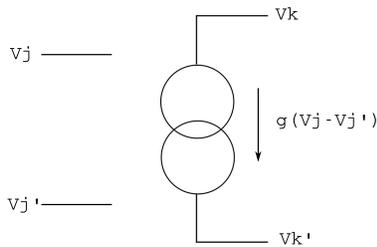}
\caption{Voltage controlled current source}
\label{fig:vccs}
\end{figure}

As long as the network consists of resistances only, the Laplacian remains symmetric, 
since for a resistive element of a conductance $c$ connected between the nodes $i$ and $j$, 
the nodal conductance matrix ${\bf L}_{R}$ can be expressed as follows:

\begin{equation}
{\bf L}_{R}=
\bordermatrix{~ & i & j \cr
i & c & -c \cr
j & -c & c \cr}.
\end{equation}
If this matrix is added to the Laplacian of the network, the symmetry of the Laplacian is maintained. 
On the other hand, if a VCCS is present in the network, the Laplacian can become non-symmetrical in cases 
when the input terminals are connected to different nodes as the output ones. In general, the nodal conductance matrix ${\bf L}_{VCCS}$ of a VCCS, 
as depicted in~\fref{fig:vccs}, can be written as

\begin{equation}
{\bf L}_{VCCS}=
\bordermatrix{~ & j & j' & k & k' \cr
j & 0 & 0 & 0 & 0 \cr
j' & 0 & 0 & 0 & 0 \cr
k & g & -g & 0 & 0 \cr
k' & -g & g & 0 & 0 \cr}.
\end{equation}

Since most of the real life circuits contain active elements that are represented by VCCSs, 
an extension of the approach published in~\cite{Wu2004} to non-symmetrical Laplacians is of great importance.

\section{The two-point resistance theorem}
\label{sec:theoretical}

%

In various applications based on matrix theory, a diagonalization of a square matrix is needed.
A real $N \times N$ matrix {\bfseries A} is diagonalizable, if there exists an invertible matrix {\bfseries P} such that ${\bf P}^{-1}{\bf AP=D}$,
where {\bfseries D} is a diagonal matrix. Here, column vectors of {\bfseries P} are right-hand eigenvectors of {\bfseries A} and row
vectors of ${\bf P}^{-1}$ are left-hand eigenvectors. 

When the matrix {\bfseries A} is symmetric, the eigenvectors can be chosen
to form an orthonormal basis of $\mathbb{R}^{N}$. In this case, {\bfseries P} is an orthogonal matrix and ${\bf P}^{-1}$ equals the transpose of
{\bfseries P}. The diagonal elements of {\bfseries D} are eigenvalues of both {\bfseries A} as {\bfseries D}.

When a diagonalizable matrix {\bfseries A} is non-symmetric, the orthogonality of eigenvectors fails.
Nevertheless, since neither the orthogonality of right eigenvectors nor that of left ones is necessary in our approach, this
problem can be overpassed. What we need, is the \textquotedblleft mutual orthogonality\textquotedblright{} of the left and right eigenvectors. Then, if $\bf \Phi$ and $\bf \Psi$ denote matrices of left and right-hand
eigenvectors of a square matrix {\bfseries A}, the matrix $\bf \Phi A \Psi$ is diagonal.
Its diagonal elements are not necessarilly the eigenvalues of {\bfseries A}.

We state our main result as a theorem. In its proof, we use a similar method as in the proof of Theorem in~\cite{Wu2004}.

From Kirchhoff's law we know that at least one eigenvalue of Laplacian matrix {\bf L} is zero. We consider this eigenvalue
simple and put $\lambda_{1}=0$. The associated left and right-hand eigenvectors are $\bm{\varphi}_{1}=(1,1,\dots,1)=\bm{\psi}^{*}_{1}$.
The condition of orthogonality of left and right-hand eigenvectors
$(\bm{\varphi}^{*}_{i}\bm{\psi}_{j}=0)\Leftrightarrow(i \neq j)$ can be fulfilled for arbitrary diagonalizable matrix that satisfies any other condition of the following theorem.
\\

\noindent{\bfseries Theorem.} {\it Consider a resistance network with $N$ nodes defined by a diagonalizable Laplacian matrix ${\bf L}$ not
necessarily symmetric. Suppose that ${\bf L}$ has $N$ real eigenvalues $\lambda_{i}$ (including the multiplicity) with left
and right-hand eigenvectors $\bm{\varphi}^{*}_{i}=(\varphi_{i1},\dots,\varphi_{iN})$ and $\bm{\psi}_{i}=(\psi_{i1},\dots,\psi_{iN})^{*}$, 
where $\bm{\varphi}_{i}$ and $\bm{\psi}_{j}$ are orthogonal if and only if $i\neq j$. Denote $\lambda_{1}$ the zero eigenvalue. 
Then the resistance between nodes $\alpha$ and $\beta$ is given by
}
\begin{equation*}
R_{\alpha\beta}=\sum^{N}_{i=2}{{1 \over \lambda_{i}\bm{\varphi}^{*}_{i}\bm{\psi}_{i}}(\varphi_{i\alpha}-\varphi_{i\beta})(\psi_{i\alpha}-\psi_{i\beta})}.
\label{eq:Rab}
\end{equation*}

\noindent{\bfseries Proof.} Since one of eigenvalues of {\bfseries L} is zero, introducing of $\epsilon{\bf L}$ to the Laplacian {\bfseries L} will be helpful.
Thus we will manipulate the Laplacian ${\bf L}(\epsilon)={\bf L}+\epsilon{\bf L}$. The eigenvalues of ${\bf L}(\epsilon)$ are $\epsilon,\lambda_{2}+\epsilon,\dots,\lambda_{N}+\epsilon$
and for a suitable value of $\epsilon$, no one of them equals zero. The left and right-hand eigenvectors are
those of {\bfseries L}. Similarly as in~\cite{Wu2004}, we will compute the inverse function ${\bm G}(\epsilon)$ of ${\bf L}(\epsilon)$.

Naturally, ${\bf L}(\epsilon)$ is not symmetric, because {\bfseries L} is not. Non-symmetric matrices {\bfseries L} and ${\bf L}(\epsilon)$ are diagonalizable as follows:
\begin{equation}
{\bf\Phi L\Psi}={\bf\Lambda}\qquad\mbox{and}\qquad{\bf\Phi L}(\epsilon)\bm{\Psi}=\bm{\Lambda}(\epsilon),
\label{eq:1_6}
\end{equation}
where ${\bf\Phi}$ is the matrix of left-hand eigenvectors (rows) and ${\bf\Psi}$ that of right-hand eigenvectors (columns) of {\bfseries L} and ${\bf L}(\epsilon)$. 
Finally, ${\bf\Lambda}(\epsilon)$ is a diagonal matrix whose diagonal elements are $(\lambda_{i}+\epsilon)\bm{\varphi}^{*}_{i}\bm{\psi}_{j}$.

Both ${\bf\Phi}$ and ${\bf\Psi}$ are invertible, because $\{\bm{\varphi}_{1},\dots,\bm{\varphi}_{N}\}$ as well as $\{\bm{\psi}_{1},\dots,\bm{\psi}_{N}\}$ define a basis of $\mathbb{R}^{N}$.
From~\eref{eq:1_6} we obtain

\begin{equation}
{\bf\Psi}^{-1}{\bf L}^{-1}(\epsilon){\bf\Phi}^{-1}={\bf\Lambda}^{-1}(\epsilon).
\label{eq:1_7}
\end{equation}
Hence

\begin{equation}
{\bf G}(\epsilon)={\bf\Psi\Lambda}^{-1}(\epsilon){\bf\Phi}.
\label{eq:1_8}
\end{equation}
Really, ${\bf G}(\epsilon)$ is an inverse matrix of ${\bf L}(\epsilon)$, because
\begin{center}
${\bf G}(\epsilon){\bf L}(\epsilon)=({\bf\Psi\Lambda}^{-1}(\epsilon){\bf\Phi})({\bf\Phi}^{-1}{\bf\Lambda}(\epsilon){\bf\Psi}^{-1})={\bf I}$,
\end{center}
as well as
\begin{center}
${\bf L}(\epsilon){\bf G}(\epsilon)=({\bf\Phi}^{-1}{\bf\Lambda}(\epsilon){\bf\Psi}^{-1})({\bf\Psi\Lambda}^{-1}(\epsilon){\bf\Phi})={\bf I}$.
\end{center}
From~\eref{eq:1_8} it follows that elements of the matrix ${\bf G}(\epsilon)$ are

\begin{equation}
{\bm G}_{\alpha\beta}(\epsilon)=\sum_{i=1}^{N}{\psi_{i\alpha}\left({1 \over (\lambda_{i}+\epsilon)\bm{\varphi}^{*}_{i}\bm{\psi}_{i}}\right)\varphi_{i\beta}={1 \over N\epsilon}+\sum_{i=2}^{N}{\psi_{i\alpha}\varphi_{i\beta}\over(\lambda_{i}+\epsilon)\bm{\varphi}^{*}_{i}\bm{\psi}_{i}}}.
\label{eq:1_9}
\end{equation}
So, the function ${\bf G}(\epsilon)$ is explicitly expressed.

Now, return to the equation ${\bf L}(\epsilon)\vec{V}(\epsilon)=\vec{I}$ and multiply it from the left by ${\bf G}(\epsilon)$. 
We obtain $\vec{V}(\epsilon)={\bf G}(\epsilon)\vec{I}$. Coordinates of the vector $\vec{V}(\epsilon)$ are

\begin{equation}
V_{i}(\epsilon)=\sum_{j=1}^{N}{G_{ij}(\epsilon)I_{j}}\quad\mbox{for}\quad i=1,2,\dots,N.
\label{eq:1_10}
\end{equation}
Denote in~\eref{eq:1_9}

\begin{equation}
g_{\alpha\beta}(\epsilon)=\sum_{i=2}^{N}{\psi_{i\alpha}\varphi_{i\beta}\over(\lambda_{i}+\epsilon)\bm{\varphi}^{*}_{i}\bm{\psi}_{i}}.
\label{eq:1_12}
\end{equation}
The substitution of $G_{ij}(\epsilon)={(N\epsilon)^{-1}}+g_{ij}(\epsilon)$ in~\eref{eq:1_10} yields
\begin{equation}
V_{i}(\epsilon)=\sum_{j=1}^{N}{\left({1\over N\epsilon}+g_{ij}(\epsilon)\right)I_{j}}={1\over N\epsilon}\sum_{j=1}^{N}{I_{j}}+\sum_{j=1}^{N}{g_{ij}(\epsilon)I_{j}}.
\label{eq:1_13}
\end{equation}
Since $\sum\limits_{i=1}^{N}{I_{j}}=0$, we have $V_{i}(\epsilon)=\sum\limits_{i=1}^{N}{g_{ij}(\epsilon)I_{j}}$ and we can put $\epsilon=0$.This gives

\begin{equation}
V_{i}=\sum_{i=1}^{N}{g_{ij}I_{j}},
\label{eq:1_14}
\end{equation}
where $V_{i}$ and $g_{ij}$ replace $V_{i}(0)$ and $g_{ij}(0)$.

\noindent We use~\eref{eq:1_14} in $R_{\alpha\beta}=I^{-1}(V_{\alpha}-V_{\beta})$ which is~\eref{eq:Wu-6}

\begin{equation}
R_{\alpha\beta}={1\over I}\left(\sum_{j=1}^{N}{g_{\alpha j}I_{j}}-\sum_{j=1}^{N}{g_{\beta j}I_{j}}\right)=\sum_{j=1}^{N}{(g_{\alpha j}-g_{\beta j}){I_{j}\over I}}.
\label{eq:1_15}
\end{equation}
Due to~\eref{eq:Wu-7} we obtain

\begin{equation}
R_{\alpha\beta}=\sum_{j=1}^{N}{(g_{\alpha j}-g_{\beta j})(\delta_{j\alpha}-\delta_{j\beta})}.
\label{eq:1_16}
\end{equation}
The last expression can be simplified as

\begin{equation}
R_{\alpha\beta}=(g_{\alpha\alpha}-g_{\beta\alpha})(1-0)+(g_{\alpha\beta}-g_{\beta\beta})(0-1).
\label{eq:1_17}
\end{equation}
We introduce~\eref{eq:1_12} with $\epsilon=0$ and~\eref{eq:1_17} becomes

\begin{equation}
R_{\alpha\beta}=\sum_{i=2}^{N}{\psi_{i\alpha}\varphi_{i\alpha}\over\lambda_{i}\bm{\varphi}_{i}^{*}\bm{\psi_{i}}}-\sum_{i=2}^{N}{\psi_{i\beta}\varphi_{i\alpha}\over\lambda_{i}\bm{\varphi}_{i}^{*}\bm{\psi_{i}}}+\sum_{i=2}^{N}{\psi_{i\beta}\varphi_{i\beta}\over\lambda_{i}\bm{\varphi}_{i}^{*}\bm{\psi_{i}}}-\sum_{i=2}^{N}{\psi_{i\alpha}\varphi_{i\beta}\over\lambda_{i}\bm{\varphi}_{i}^{*}\bm{\psi_{i}}}.
\end{equation}
A simple rearrangement of the last relation returns~\eref{eq:novy_vv}.\qed

\section{Examples}
\label{sec:examples}

In the following examples, the proposed resistance calculation approach is demonstrated on simple networks (circuits). 
The calculations were conducted using {\it GNU Octave} language. The circuit analysis results were obtained from {\it Cadence IC} design environment.
\\

\noindent{\bfseries Example 1 (simple eigenvalues).} 

\begin{figure}[h!]
\includegraphics[width=5cm]{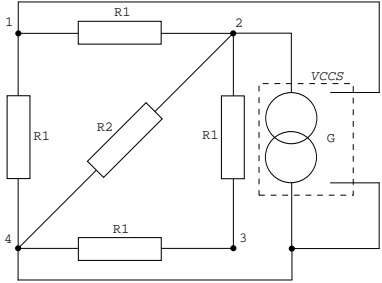}
\caption{A network of 4 nodes with VCCS}
\label{fig:bridge}
\end{figure}

Consider the network depicted in~\fref{fig:bridge}. For this network, the nodal conductance matrix,
which is the Laplacian {\bf L}, is the following non-symmetric matrix, where $c_{1}=1/R_{1}\mbox{ and }c_{2}=1/R_{2}$:

$  \mathbf{L}=\begin{blockarray}{cccc}
  & & & & \\
    \begin{block}{(rrrr)}
2c_{1} & -c_{1} & 0 & -c_{1} \\
-c_{1}+G & 2c_{1}+c_{2} & -c_{1} & -c_{2}-G \\
0 & -c_{1} & 2c_{1} & -c_{1} \\
-c_{1}-G & -c_{2} & -c_{1} & 2c_{1}+c_{2}+G \\
    \end{block}
  \end{blockarray}$

The eigenvalues are $\lambda_{1}=0, \lambda_{2}=2c_{1}, \lambda_{3}=4c_{1}, \lambda_{4}=2c_{1}+2c_{2}+G$.
Since the eigenvectors are not required to have the norm 1, we take those with the least complicated entries.
The left and right-hand eigenvectors are

$  \mathbf{\Phi}^{T}=\begin{blockarray}{cccc}
    \matindex{$\bm{\varphi}^{*}_{1}$} & \matindex{$\bm{\varphi}^{*}_{2}$} & \matindex{$\bm{\varphi}^{*}_{3}$} & \matindex{$\bm{\varphi}^{*}_{4}$} \\
    \begin{block}{(rrrr)}
      1 & 1 & 1 & G(c_{1}-2c_{2}-G) \\
      1 & 0 & -1 & (c_{1}-c_{2})(2c_{2}+G) \\
      1 & -1 & 1 & -c_{1}G \\
      1 & 0 & -1 & (2c_{2}+G)(-c_{1}+c_{2}+G) \\
    \end{block}
  \end{blockarray}$

$  \mathbf{\Psi}=\begin{blockarray}{cccc}
    \matindex{$\bm{\psi}_{1}$} & \matindex{$\bm{\psi}_{2}$} & \matindex{$\bm{\psi}_{3}$} & \matindex{$\bm{\psi}_{4}$} \\
    \begin{block}{(rrrr)}
      1 & -2c_{2}-G & -2c_{1}+2c_{2}+G & 0 \\
      1 & G & 2c_{1}-2c_{2}-3G & 1 \\
      1 & 2c_{2}+G & -2c_{1}+2c_{2}+G & 0 \\
      1 & -G & 2c_{1}-2c_{2}+G & -1 \\
    \end{block}
  \end{blockarray}$





%
Using~\eref{eq:novy_vv}, the resistance between two arbitrary nodes of the network can now be calculated. 
Below, examples of resistance calculation follow.
\begin{eqnarray*}
R_{13}=&
{{1 \over \lambda_{2}\bm{\varphi}^{*}_{2}\bm{\psi}_{2}}(\varphi_{21}-\varphi_{23})(\psi_{21}-\psi_{23})}+
{{1 \over \lambda_{3}\bm{\varphi}^{*}_{3}\bm{\psi}_{3}}(\varphi_{31}-\varphi_{33})(\psi_{31}-\psi_{33})}+\\
&+{{1 \over \lambda_{4}\bm{\varphi}^{*}_{4}\bm{\psi}_{4}}(\varphi_{41}-\varphi_{43})(\psi_{41}-\psi_{43})}=R_{1},
\end{eqnarray*}
\begin{eqnarray*}
R_{24}=&
{{1 \over \lambda_{2}\bm{\varphi}^{*}_{2}\bm{\psi}_{2}}(\varphi_{22}-\varphi_{24})(\psi_{22}-\psi_{24})}+
{{1 \over \lambda_{3}\bm{\varphi}^{*}_{3}\bm{\psi}_{3}}(\varphi_{32}-\varphi_{34})(\psi_{32}-\psi_{34})}+\\
&+{{1 \over \lambda_{4}\bm{\varphi}^{*}_{4}\bm{\psi}_{4}}(\varphi_{42}-\varphi_{44})(\psi_{42}-\psi_{44})}
={2R_{1}R_{2}\over 2R_{1}+2R_{2}+GR_{1}R_{2}}.
\end{eqnarray*}
Set the values of parameters $R_{1}=200\Omega, R_{2}=4000\Omega, G=0.03S$. 
All the possible two-node resistances constitute an $N \times N$ upper triangular matrix where $N$ is the network node count.
The matrices of calculated resistance $\mathbf{R}_{c}$ and resistance from circuit analysis $\mathbf{R}_{a}$ are shown below.
One can observe that $\mathbf{R}_{c}=\mathbf{R}_{a}$.

$  \mathbf{R}_{c}=\begin{blockarray}{cccc}
  & & & & \\
    \begin{block}{(rrrr)}
0 & 1.864E+02 & 2.000E+02 & 3.827E+01\\
0 & 0 & 1.124E+02 & 4.938E+01\\
0 & 0 & 0 & 1.124E+02\\
0 & 0 & 0 & 0\\
    \end{block}
  \end{blockarray}
$

$  \mathbf{R}_{a}=\begin{blockarray}{cccc}
  & & & & \\
    \begin{block}{(rrrr)}
0 & 1.864E+02 & 2.000E+02 & 3.827E+01\\
0 & 0 & 1.124E+02 & 4.938E+01\\
0 & 0 & 0 & 1.124E+02\\
0 & 0 & 0 & 0\\
    \end{block}
  \end{blockarray}$

\noindent{\bfseries Example 2 (multiple eigenvalues).} 

\begin{figure}[h!]
\includegraphics[width=5cm]{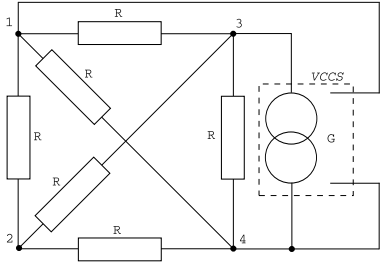}
\caption{A network of 4 nodes with VCCS}
\label{fig:multibridge}
\end{figure}

Consider the network depicted in~\fref{fig:multibridge}. Let $c_{1}=1/R$.
For this network, the following non-symmetrical nodal conductance matrix, which is the Laplacian $\mathbf{L}$, can be created:

\noindent$
  \mathbf{L}=\begin{blockarray}{cccc}
  & & & & \\
    \begin{block}{(rrrr)}
3c_{1} & -c_{1} & -c_{1} & -c_{1} \\
-c_{1} & 3c_{1} & -c_{1} & -c_{1} \\
-c_{1}+G & -c_{1} & 3c_{1} & -c_{1}-G \\
-c_{1}-G & -c_{1} & -c_{1} & 3c_{1}+G \\
    \end{block}
  \end{blockarray}
$\\
The eigenvalues are $\lambda_{1}=0,\lambda_{2}=4c_{1},\lambda_{3}=4c_{1},\lambda_{4}=4c_{1}+G$. Note the multiple eigenvalue $\lambda_{2}=\lambda_{3}=4c_{1}$. 
Left and right-hand eigenvectors $\bm{\varphi}_{i}\mbox{ and }\bm{\psi}_{i}$ of this matrix are as follows:

\noindent$  \mathbf{\Phi}^{T}=\begin{blockarray}{cccc}
    \matindex{$\bm{\varphi}^{*}_{1}$} & \matindex{$\bm{\varphi}^{*}_{2}$} & \matindex{$\bm{\varphi}^{*}_{3}$} & \matindex{$\bm{\varphi}^{*}_{4}$} \\
    \begin{block}{(rrrr)}
      1 & -1 & 1 & 1 \\
      1 & -1 & -1 & 0 \\
      1 & 1 & 0 & 0 \\
      1 & 1 & 0 & -1 \\
    \end{block}
  \end{blockarray}
\quad  \mathbf{\Psi}=\begin{blockarray}{cccc}
    \matindex{$\bm{\psi}_{1}$} & \matindex{$\bm{\psi}_{2}$} & \matindex{$\bm{\psi}_{3}$} & \matindex{$\bm{\psi}_{4}$} \\
    \begin{block}{(rrrr)}
      1 & 1 & 1 & 0 \\
      1 & 1 & -1 & 0 \\
      1 & -3 & -1 & 1 \\
      1 & 1 & 1 & -1 \\
    \end{block}
  \end{blockarray}$\\
As can be observed, the left and right-hand eigenvectors are independent of $c_{1} \mbox{ and }G$.
An example calculation of the first row of the resistance matrix $\mathbf{R}_{c}$ is given below.
We can use substitution $c_{1}=1/R$ for this calculation.
\begin{equation*}
{1 \over \lambda_{2}\bm{\varphi}^{*}_{2}\bm{\psi}_{2}}=-{R \over 16},\qquad
{1 \over \lambda_{3}\bm{\varphi}^{*}_{3}\bm{\psi}_{3}}={R \over 8},\qquad
{1 \over \lambda_{4}\bm{\varphi}^{*}_{4}\bm{\psi}_{4}}={R \over 4+GR},
\end{equation*}

\begin{eqnarray*}
R_{12} &=-{R \over 16}(-1+1)(1-1)+{R \over 8}(1+1)(1+1)+{R \over 4+GR}(1-0)(0-0)\\
&=0+{R \over 8}4+0={R \over 2},
\end{eqnarray*}

\begin{eqnarray*}
R_{13}&=-{R \over 16}(-1-1)(1+3)+{R \over 8}(1-0)(1+1)+{R \over 4+GR}(1-0)(0-1)\\
&={R \over 2}+{R \over 4}-{R \over 4+GR}=R\left({3 \over 4}-{1\over 4+GR}\right),
\end{eqnarray*}

\begin{eqnarray*}
R_{14}&=-{R \over 16}(-1-1)(1-1)+{R \over 8}(1-0)(1-1)+{R \over 4+GR}(1+1)(0+1)\\
&=0+0+{2R \over 4+GR}.
\end{eqnarray*}

Now, set the values of parameters $R=1\Omega$ and $G=2S$. The  values of resistances are
$R_{12}={1/2}, R_{13}={7/12} \mbox{ and }R_{14}={1/3}$.
The complete resulting matrices of calculated resistance $\mathbf{R}_{c}$ and resistance from circuit analysis $\mathbf{R}_{a}$ are shown below.
One can observe that $\mathbf{R}_{c}=\mathbf{R}_{a}$.

\noindent$  \mathbf{R}_{c}=\begin{blockarray}{cccc}
  & & & \\
    \begin{block}{(rrrr)}
0 & 1/2 & 7/12 & 1/3 \\
0 & 0 & 1/2 & 5/12 \\
0 & 0 & 0 & 1/3 \\
0 & 0 & 0 & 0 \\
    \end{block}
  \end{blockarray}
\quad \mathbf{R}_{a}=\begin{blockarray}{cccc}
  & & & \\
    \begin{block}{(rrrr)}
0 & 5.000E-01 & 5.833E-01 & 3.333E-01 \\
0 & 0 & 5.000E-01 & 4.167E-01 \\
0 & 0 & 0 & 3.333E-01 \\
0 & 0 & 0 & 0 \\
    \end{block}
  \end{blockarray}$

\noindent{\bfseries Example 3 (operational amplifier).} 

\begin{figure}[h!]
\includegraphics[width=7cm]{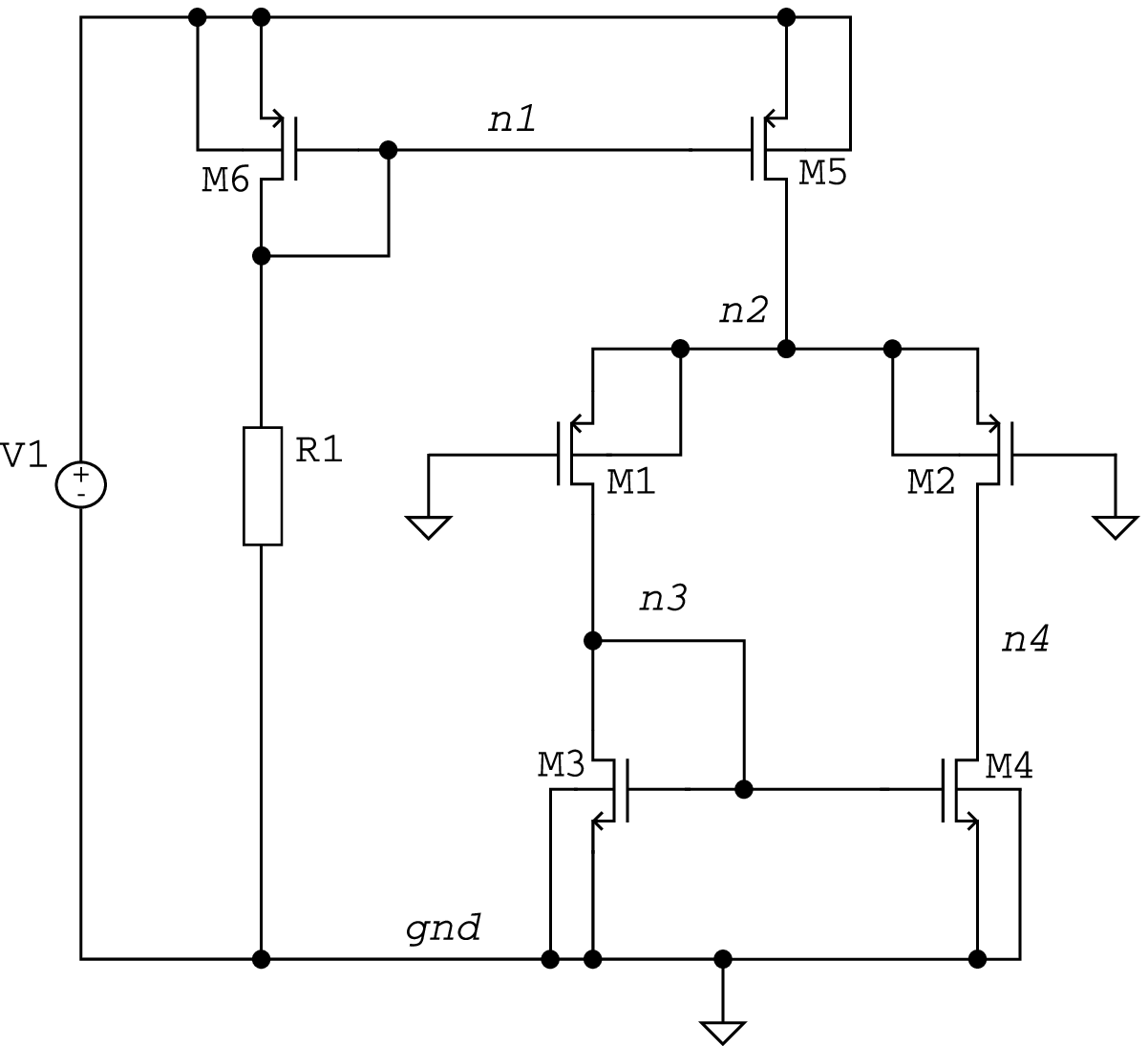}
\caption{Circuit schematic of an operational amplifier}
\label{fig:opamp}
\end{figure}

The circuit depicted in~\fref{fig:opamp} represents a typical topology of an operational amplifier 
that is one of the most frequently used analog building block.

This particular circuit has been designed for a standard $0.35\mu m$ mixed-signal CMOS process. 
Parameters of the devices obtained from a typical $.op$ (operating point) analysis are shown in \tref{tab:opamp}.

\begin{table}[h]
\caption{Device parameters of the operational amplifier}
\label{tab:opamp}
\begin{indented}
\lineup
\item[]\begin{tabular}{@{}*{4}{l}}
\br                              
&gm[S]&gds[S]&Ron[$\Omega$]\cr 
\mr
M1&$5.155032E-04$&$2.816268E-03$&-\cr
M2&$5.155032E-04$&$2.816268E-03$&-\cr
M3&$5.953932E-04$&$3.403570E-06$&-\cr
M4&$5.953932E-04$&$3.403570E-06$&-\cr
M5&$4.195651E-03$&$9.344114E-06$&-\cr
M6&$5.247305E-04$&$2.065995E-06$&-\cr
R1&-&-&$6.847545E+04$\cr
\br
\end{tabular}
\end{indented}
\end{table}

\begin{figure}[h!]
\includegraphics[width=4cm]{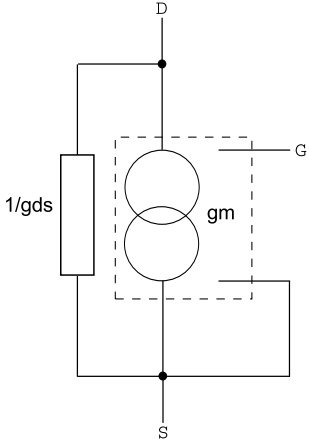}
\caption{Simplified electrical model of a MOS transistor}
\label{fig:mos_model}
\end{figure}

These parameters are used to create an electrical model of the circuit, where only resistances and voltage-controlled current sources (VCCS) are used. 
MOS transistor is modeled by a VCCS, which is represented by a transconductance {\it gm} and a resistor of value {\it 1/gds} connected between its output terminals, as depicted in~\fref{fig:mos_model}.

Using the parameters from \tref{tab:opamp}, 
the following non-symmetrical nodal admittance matrix ${\mathbf L}$ of this circuit can be created for the internal nodes:

\noindent$  \mathbf{L}=\begin{blockarray}{ccccc}
  & & & & \\
    \begin{block}{(rrrrr)}
5.414E-04 & 0.000E+00 & 0.000E+00 & 0.000E+00 & -5.414E-04 \\
4.196E-03 & 6.673E-03 & -2.816E-03 & -2.816E-03 & -5.236E-03 \\
0.000E+00 & -3.332E-03 & 3.415E-03 & 0.000E+00 & -8.329E-05 \\
0.000E+00 & -3.332E-03 & 5.954E-04 & 2.820E-03 & -8.329E-05 \\
-4.737E-03 & -9.344E-06 & -1.194E-03 & -3.404E-06 & 5.944E-03 \\
    \end{block}
  \end{blockarray}
$\\
The node $Vdd$ is excluded from this matrix, since when assuming the impedance between two nodes, 
all the voltage and current sources have to be substituted by their equivalent circuits. 
Since no parasitic resistance is taken into account, a voltage source is substituted by a short, while a current source is represented by an open. 
In our case, it means that the node $Vdd$ becomes shorted to node $gnd$. The order of network nodes assumed for this matrix is $n1, n2, n3, n4, gnd$.

For the nodal admittance matrix, which is the Laplacian $\mathbf{L}$, the eigenvalues obtained are $\lambda_{1}=0,\lambda_{2}=1.617872E-04,
\lambda_{3}=2.819672E-03,\lambda_{4}=8.008586E-03, \lambda_{5}=8.402968E-03$.

The corresponding left and right-hand eigenvectors (in columns) of the matrix $\mathbf{L}$ are as follows:

{$\mathbf{\Phi}^{T}=\begin{blockarray}{ccccc}
    \matindex{$\bm{\varphi}^{*}_{1}$} & \matindex{$\bm{\varphi}^{*}_{2}$} & \matindex{$\bm{\varphi}^{*}_{3}$} & \matindex{$\bm{\varphi}^{*}_{4}$} & \matindex{$\bm{\varphi}^{*}_{5}$} \\
    \begin{block}{(rrrrr)}
 0.447214 & 0.890916 & 0 & 0.606347 & -0.596298 \\
 0.447214 & -0.250715 & 0 & 0.233943 & -0.286327 \\
 0.447214 & -0.223688 & -0.707107 & 0.034793 & 0.002639 \\
 0.447214 & -0.265848 & 0.707107 & -0.126481 & 0.143977 \\
 0.447214 & -0.150665 & 0 & -0.748603 & 0.736009 \\
    \end{block}
  \end{blockarray}$

$\mathbf{\Psi}=\begin{blockarray}{ccccc}
    \matindex{$\bm{\psi}_{1}$} & \matindex{$\bm{\psi}_{2}$} & \matindex{$\bm{\psi}_{3}$} & \matindex{$\bm{\psi}_{4}$} & \matindex{$\bm{\psi}_{5}$} \\
    \begin{block}{(rrrrr)}
 0.447214 & -0.638056 & -0.038249 & 0.023987 & -0.017999 \\
 0.447214 & 0.363601 & 0.098610 & -0.654286 & 0.698457 \\
 0.447214 & 0.360920 & 0.574328 & 0.480565 & -0.470914 \\
 0.447214 & 0.360920 & -0.795646 & 0.480565 & -0.470914 \\
 0.447214 & -0.447385 & 0.160957 & -0.330831 & 0.261369 \\
    \end{block}
  \end{blockarray}$
}\\
The resulting matrices of the numerically calculated resistance $\mathbf{R}_{c}$ and the resistance from the circuit analysis $\mathbf{R}_{a}$ are shown below.
One can observe that $\mathbf{R}_{c}=\mathbf{R}_{a}$.

{$\mathbf{R}_{c}=\begin{blockarray}{ccccc}
  & & & & \\
    \begin{block}{(rrrrr)}
0 & 9.276302E+03 & 9.098426E+03 & 9.448754E+03 & 1.847062E+03 \\
0 & 0 & 3.000575E+02 & 3.546762E+02 & 8.490922E+02 \\
0 & 0 & 0 & 7.093024E+02 & 8.317060E+02 \\
0 & 0 & 0 & 0 & 1.182034E+03 \\
0 & 0 & 0 & 0 & 0 \\
    \end{block}
  \end{blockarray}$

$\mathbf{R}_{a}=\begin{blockarray}{ccccc}
  & & & & \\
    \begin{block}{(rrrrr)}
0 & 9.276302E+03 & 9.098426E+03 & 9.448754E+03 & 1.847062E+03 \\
0 & 0 & 3.000575E+02 & 3.546762E+02 & 8.490922E+02 \\
0 & 0 & 0 & 7.093024E+02 & 8.317060E+02 \\
0 & 0 & 0 & 0 & 1.182034E+03 \\
0 & 0 & 0 & 0 & 0 \\
    \end{block}
  \end{blockarray}$
}
\newpage
\section{Conclusion}

A formula for resistance computation between two arbitrary nodes of a network with non-symmetric Laplacian matrix {\bf L} has been formulated and proved in this paper.
The calculation is based on eigenvalues and eigenvectors of the Laplacian matrix.
The application of the formula is demonstrated on three examples of networks (circuits): two resistive bridges and an operational amplifier.
The calculated results are compared to results from a circuit simulator.

In electric circuit theory, resistor networks are a special case of impedance networks. 
In general, these are defined by a complex Laplacian matrix {\bf L}. 
Our numerical results indicate, that Theorem formulated in this paper holds also for systems with complex Laplacian matrix.
Therefore,~\eref{eq:novy_vv} should be applicable to computation of the complex impedance as well.

%

\end{document}